# О связи моделей дискретного выбора с разномасштабными по времени популяционными играми загрузок

*Гасников А.В., Гасникова Е.В., Мациевский С.В., Усик И.В.*
(Центр исследований транспортной политики, Институт экономики транспорта и транспортной политики, НИУ ВШЭ; Балтийский федеральный университет им. И. Канта; ПреМоЛаб ФУПМ МФТИ)

**Аннотация**
В работе предложен универсальный прямо-двойственный способ описания равновесия в иерархических популяционных играх загрузок. В основе подхода лежит иерархия вложенных друг в друга транспортных сетей и соответствующие этим сетям разномасштабные (по времени) логит динамики, отражающие ограниченную рациональность агентов. Поиск равновесия сводится к решению мгногоуровневой задачи выпуклой оптимизации. Описанные в статье результаты могут быть использованы при описании и численном поиске равновесий (стохастических равновесий) во всех известных многостадийных моделях транспортных потоков.

## 1. Введение

В работах [1 – 4] было анонсировано, что в последующем цикле публикаций будет приведен общий способ вариационного описания равновесий (стохастических равновесий) в популярных моделях распределения транспортных потоков. Также отмечалось, что планируется предложить эффективные численные методы поиска таких равновесий. В данной работе предпринята попытка погрузить известные нам подходы к многостадийному моделированию потоков на иерархических сетях (реальных транспортных сетях или сетях принятия решение – не важно) в одну общую схему, сводящую поиск равновесия к решению многоуровневой задачи выпуклой оптимизации. В основе схемы получения вариационного принципа для описания равновесия лежит популяционная игра загрузки с соответствующими логит динамиками (отвечающими моделям дискретного выбора [5]) пользователей на каждом уровне иерархии [6] (см. п. 2). Для решения описанной задачи выпуклой оптимизации в статье изучается двойственная задача, представляющая самостоятельный интерес (см. п. 3). Основным инструментом изучения двойственной задачи является аппарат характеристических функций на графе [1, 7, 8] и ускоренные прямо-двойственные методы в композитном варианте [9, 10].

Отметим, что общность результатов статьи достигается за счет введения большого числа параметров, которые можно вырождать, стремя их к нулю или бесконечности. Игра на выборе этих параметров позволяет, например, получать различные многостадийные модели транспортных потоков [2 – 4, 11, 12]. Приводимые далее результаты можно обобщать и на потоки товаров, в случае, когда имеется более одного наименования товара [12]. Однако в данной статье мы не будем касаться этого обобщения. Мы также не планируем приводить конкретные примеры получения многостадийных транспортных моделей согласно изложенной в статье общей схеме. Всему этому планируется посвятить отдельные публикации.

## 2. Постановка задачи

Рассмотрим транспортную сеть, заданную ориентированным графом $\Gamma^1 = \left\langle V^1, E^1 \right\rangle$. Часть его вершин $O^1 \subseteq V^1$ является источниками, часть стоками $D^1 \subseteq V^1$. Множество пар источник-сток, обозначим $OD^1 \subseteq O^1 \otimes D^1$. Пусть каждой паре $w^1 \in OD^1$ соответствует своя корреспонденция: $d^1_{w^1} := d^1_{w^1} \cdot M$ ($M \gg 1$) пользователей, которые хотят в единицу времени перемещаться из источника в сток, соответствующих заданной корреспонденции $w^1$. Пусть ребра $\Gamma^1$ разделены на два типа $E^1 = \tilde{E}^1 \coprod \bar{E}^1$. Ребра типа $\tilde{E}^1$ характеризуются неубывающими функциями затрат $\tau^1_{e^1}\left(f^1_{e^1}\right) := \tau^1_{e^1}\left(f^1_{e^1}/M\right)$. Затраты $\tau^1_{e^1}\left(f^1_{e^1}\right)$ несут те пользователи, которые используют в своем пути ребро $e^1 \in \tilde{E}^1$, в предположении, что поток пользователей по этому ребру равен $f^1_{e^1}$. Пары вершин, задающие ребра типа $\bar{E}^1$, являются, в свою очередь, парами источник-сток $OD^2$ (с корреспонденциями $d^2_{w^2} = f^1_{e^1}$, $w^2 = e^1 \in \bar{E}_1$) в транспортной сети следующего уровня, $\Gamma^2 = \left\langle V^2, E^2 \right\rangle$, ребра которой, в свою очередь, разделены на два типа $E^2 = \tilde{E}^2 \coprod \bar{E}^2$. Ребра типа $\tilde{E}^2$ характеризуются неубывающими функциями затрат $\tau^2_{e^2}\left(f^2_{e^2}\right) := \tau^2_{e^2}\left(f^2_{e^2}/M\right)$. Затраты $\tau^2_{e^2}\left(f^2_{e^2}\right)$ несут те пользователи, которые используют в своем пути ребро $e^2 \in \tilde{E}^2$, в предположении, что поток пользователей по этому ребру равен $f^2_{e^2}$. Пары вершин, задающие ребра типа $\bar{E}^2$, являются, в свою очередь, парами источник-сток $OD^3$ (с корреспонденциями $d^3_{w^3} = f^2_{e^2}$, $w^3 = e^2 \in \bar{E}^2$) в транспортной сети более высокого уровня, $\Gamma^3 = \left\langle V^3, E^3 \right\rangle$, … и т.д. Будем считать, что всего имеется $m$ уровней: $\tilde{E}^m = E^m$. Обычно в приложениях число $m$ небольшое [1 – 5]: 2 – 10.

Каждый пользователь в графе $\Gamma^1$ выбирает путь $p^1_{w^1} \in P^1_{w^1}$ (последовательный набор проходимых пользователем ребер), соответствующий его корреспонденции $w^1 \in OD^1$ ($P^1_{w^1}$ – множество всех путей, отвечающих в $\Gamma^1$ корреспонденции $w^1$). Задав $p^1_{w^1}$ можно однозначно восстановить ребра типа $\bar{E}^1$, входящие в этот путь. На каждом из этих ребер $w^2 \in \bar{E}^1$ пользователь может выбирать свои пути $p^2_{w^2} \in P^2_{w^2}$ ($P^2_{w^2}$ – множество всех путей, отвечающих в $\Gamma^2$ корреспонденции $w^2$), … и т.д. Пусть каждый пользователь сделал свой выбор. Обозначим через $x^1_{p^1}$ величину потока пользователей по пути $p^1 \in P^1 = \coprod_{w^1 \in OD^1} P^1_{w^1}$, $x^2_{p^2}$ – величина потока пользователей по пути $p^2 \in P^2 = \coprod_{w^2 \in OD^2} P^2_{w^2}$, … и т.д. Заметим, что

$$x^k_{p^k_{w^k}} \geq 0, \; p^k_{w^k} \in P^k_{w^k}, \; \sum_{p^k_{w^k} \in P^k_{w^k}} x^k_{p^k_{w^k}} = d^k_{w^k}, \; w^k \in OD^k, \; k = 1, ..., m,$$

что для компактности мы будем далее записывать

$$\left\{x^k_{p^k_{w^k}}\right\}_{p^k_{w^k}\in P^k_{w^k}} \in S_{\left|P^k_{w^k}\right|}\left(d_{w^k}\right).$$

Отметим, что здесь и везде в дальнейшем

$$w^{k+1}\left(=e^k\right)\in OD^{k+1}\left(=\bar{E}^k\right),\; d^{k+1}_{w^{k+1}} = f^k_{e^k},\; k=1,...,m-1.$$

Введем для графа $\Gamma^k$ и множества путей $P^k$ матрицу (Кирхгофа)

$$\Theta^k = \left\|\delta_{e^k p^k}\right\|_{e^k\in E^k,\, p^k\in P^k},\; \delta_{e^k p^k} = \begin{cases}1, & e^k\in p^k\\ 0, & e^k\notin p^k\end{cases},\; k=1,...,m.$$

Тогда вектора потоков на ребрах $f^k$ на графе $\Gamma^k$ однозначно определяются векторами потоков на путях $x^k = \left\{x^k_{p^k}\right\}_{p^k\in P^k}$:

$$f^k = \Theta^k x^k,\; k=1,...,m.$$

Обозначим через

$$x = \left\{x^k\right\}_{k=1}^m,\; f = \left\{f^k\right\}_{k=1}^m,\; \Theta = \mathrm{diag}\left\{\Theta^k\right\}_{k=1}^m,$$

$$d^k = \left\{d^k_{w^k}\right\}_{w^k\in OD^k},\; X^k\left(d^k\right) = \coprod_{w^k\in OD^k} S_{\left|P^k_{w^k}\right|}\left(d^k_{w^k}\right),\; X = \coprod_{k=1}^m X^k\left(d^k\right),$$

а через

$$\breve{p}^k_{w^k} = \left(p^k_{w^k}, \left\{p^{k+1}_{w^{k+1}}\right\}_{w^{k+1}\in p^k_{w^k}\cap\bar{E}^k}, ...\right),\; k=1,...,m$$

полное описание возможного пути (в графе $\Gamma^k$ и графах следующих уровней), соответствующего корреспонденции $w^k\in OD^k$. Множество всех таких путей будем обозначать $\breve{P}^k_{w^k}$. Введем также множество путей $\breve{P}^k = \coprod_{w^k\in OD^k}\breve{P}^k_{w^k}$ и соответствующий вектор распределения потоков по этим путям $x_{\breve{P}^k}$. Определим функции затрат пользователей на пути $\breve{p}^k_{w^k}$ по индукции:

$$G^m_{\breve{p}^m_{w^m}}\left(x_{\breve{P}^m}\right) = \sum_{e^m\in\tilde{E}^m}\delta_{e^m p^m}\tau^m_{e^m}\left(f^m_{e^m}\right),$$

$$G^k_{\breve{p}^k_{w^k}}\left(x_{\breve{P}^k}\right) = \sum_{e^k\in\tilde{E}^k}\delta_{e^k\breve{p}^k_{w^k}}\tau^k_{e^k}\left(f^k_{e^k}\right) + \sum_{w^{k+1}\in\bar{E}^k} G^{k+1}_{\breve{p}^{k+1}_{w^{k+1}}}\left(x_{\breve{P}^{k+1}}\right),\; k=1,...,m-1.$$

Опишем *марковскую логит динамику* (также говорят гиббсовскую динамику) в повторяющейся игре загрузки графа транспортной сети [6, 8]. Пусть имеется $TN$ шагов ($N\gg 1$). Каждый пользователь транспортной сети, использовавший на шаге $t$ путь $\breve{p}^1_{w^1}$, независимо от остальных, на шаге $t+1$ (все введенные новые параметры положительны)

- с вероятностью

$$\frac{\lambda^1}{N}\frac{\exp\left(-G^{t,1}_p/\gamma^1\right)}{\sum\limits_{\tilde{p}\in\breve{P}^1_{w^1}}\exp\left(-G^{t,1}_{\tilde{p}}/\gamma^1\right)}.$$

пытается изменить свой путь $\breve{p}^1_{w^1}$ на $p\in\breve{P}^1_{w^1}$, где $G^{t,1}_p = G^1_p\left(x^t_{\breve{P}^1}\right)$ – затраты на пути $p$ на шаге $t$ ($G^{0,1}_p\equiv 0$);

- равновероятно выбирает $w^2 \in p^1_{w^1} \cap \bar{E}^1$ и затем с вероятностью

$$\frac{\lambda^2 \left| p^1_{w^1} \cap \bar{E}^1 \right|}{N} \frac{\exp\left(-G^{t,2}_p / \gamma^2\right)}{\sum_{\tilde{p} \in \breve{P}^2_{w^2}} \exp\left(-G^{t,2}_{\tilde{p}} / \gamma^2\right)}$$

пытается изменить в своем пути $\breve{p}^1_{w^1} = \left( p^1_{w^1}, \left\{ \breve{p}^2_{w^2} \right\}_{w^2 \in p^1_{w^1} \cap \bar{E}^1} \right)$ участок пути $\breve{p}^2_{w^2}$, выбирая путь $p \in \breve{P}^2_{w^2}$, где $G^{t,2}_p = G^2_p\left(x^t_{\breve{P}^2}\right)$ – затраты на пути $p$ на шаге $t$ ($G^{0,2}_p \equiv 0$);

- ... и т.д.;
- с вероятностью

$$1 - \sum_{k=1}^{m} \lambda^k \Big/ N$$

решает не менять тот путь, который использовал на шаге $t$.

Такая динамика отражает ограниченную рациональность агентов (см. замечание 5 п. 3), и часто используется в теории дискретного выбора [5] и популяционной теории игр [6]. В основном нас будет интересовать поведение такой системы в предположении

$$\lambda^2 / \lambda^1 \to \infty,\ \lambda^3 / \lambda^2 \to \infty,\ \ldots,\ \lambda^m / \lambda^{m-1} \to \infty,\ N / \lambda^m \to \infty. \tag{1}$$

Эта марковская динамика в пределе $N \to \infty$ превращается в марковскую динамику в непрерывном времени [13]. Далее мы, как правило, будем считать, что такой предельный переход был осуществлен.

В пределе $M \to \infty$ эта динамика (концентраций) описывается зацепляющейся системой обыкновенных дифференциальных уравнений (*СОДУ*)

$$\frac{dx_{\breve{p}^1_{w^1}}}{dt} = \lambda^1 \cdot \left( d^1_{w^1} \frac{\exp\left(-G^1_{\breve{p}^1_{w^1}}\left(x_{\breve{P}^1}\right) / \gamma^1\right)}{\sum_{\tilde{p} \in \breve{P}^1_{w^1}} \exp\left(-G^1_{\tilde{p}}\left(x_{\breve{P}^1}\right) / \gamma^1\right)} - x_{\breve{p}^1_{w^1}} \right),\ \breve{p}^1_{w^1} \in \breve{P}^1_{w^1},\ w^1 \in OD^1,$$

$$\frac{dx_{\breve{p}^2_{w^2}}}{dt} = \lambda^2 \cdot \left( d^2_{w^2} \frac{\exp\left(-G^2_{\breve{p}^2_{w^2}}\left(x_{\breve{P}^2}\right) / \gamma^2\right)}{\sum_{\tilde{p} \in \breve{P}^2_{w^2}} \exp\left(-G^2_{\tilde{p}}\left(x_{\breve{P}^2}\right) / \gamma^2\right)} - x_{\breve{p}^2_{w^2}} \right),\ \breve{p}^2_{w^2} \in \breve{P}^2_{w^2},\ w^2 \in OD^2 = \bar{E}^1,$$

…………………………………………………………………………

Применяя по индукции (в виду условия (1)) теорему Тихонова к этой СОДУ [14, 15], можно получить описание аттрактора СОДУ – глобально устойчивой (при $T \to \infty$) неподвижной точки. Для того чтобы это сделать, введем обозначение

$$\sigma^k_{e^k}\left(f^k_{e^k}\right) = \int_0^{f^k_{e^k}} \tau^k_{e^k}(z)\,dz,\ k = 1,\ldots,m.$$

Рассмотрим задачу

$$\Psi(x, f) := \Psi^1(x) = \sum_{e^1 \in \tilde{E}^1} \sigma^1_{e^1}\left(f^1_{e^1}\right) + \Psi^2(x) + \gamma^1 \sum_{w^1 \in OD^1} \sum_{p^1 \in P^1_{w^1}} x^1_{p^1} \ln\left(x^1_{p^1} / d^1_{w^1}\right) \to \min_{f = \Theta x,\, x \in X}, \tag{2}$$

$$\Psi^2(x) = \sum_{e^2 \in \tilde{E}^2} \sigma^2_{e^2}\left(f^2_{e^2}\right) + \Psi^3(x) + \gamma^2 \sum_{w^2 \in \bar{E}^1} \sum_{p^2 \in P^2_{w^2}} x^2_{p^2} \ln\left(x^2_{p^2} / d^2_{w^2}\right),\ d^2_{w^2} = f^1_{w^2},$$

……………………………………………………………………………………..

$$\Psi^k(x) = \sum_{e^k \in \tilde{E}^k} \sigma_{e^k}^k\left(f_{e^k}^k\right) + \Psi^{k+1}(x) + \gamma^k \sum_{w^k \in \bar{E}^{k-1}} \sum_{p^k \in P_{w^k}^k} x_{p^k}^k \ln\left(x_{p^k}^k / d_{w^k}^k\right),\ d_{w^{k+1}}^{k+1} = f_{w^{k+1}}^k,$$

……………………………………………………………………………………..

$$\Psi^m(x) = \sum_{e^m \in E^m} \sigma_{e^m}^m\left(f_{e^m}^m\right) + \gamma^m \sum_{w^m \in \bar{E}^{m-1}} \sum_{p^m \in P_{w^m}^m} x_{p^m}^m \ln\left(x_{p^m}^m / d_{w^m}^m\right),\ d_{w^m}^m = f_{w^m}^{m-1}.$$

Эта задача эквивалентна следующей цепочке зацепляющихся задач выпуклой (многоуровневой [16]) оптимизации

$$\Phi^1\left(d^1\right) = \min_{\substack{f^1 = \Theta^1 x^1,\ x^1 \in X^1\left(d^1\right) \\ d_{e^1}^2 = f_{e^1}^1,\ e^1 \in \bar{E}^1}} \left\{ \sum_{e^1 \in \tilde{E}^1} \sigma_{e^1}^1\left(f_{e^1}^1\right) + \Phi^2\left(d^2\right) + \gamma^1 \sum_{w^1 \in OD^1} \sum_{p^1 \in P_{w^1}^1} x_{p^1}^1 \ln\left(x_{p^1}^1 / d_{w^1}^1\right) \right\}, \quad (3)$$

$$\Phi^2\left(d^2\right) = \min_{\substack{f^2 = \Theta^2 x^2,\ x^2 \in X^2\left(d^2\right) \\ d_{e^2}^3 = f_{e^2}^2,\ e^2 \in \bar{E}^2}} \left\{ \sum_{e^2 \in \tilde{E}^2} \sigma_{e^2}^2\left(f_{e^2}^2\right) + \Phi^3\left(d^3\right) + \gamma^2 \sum_{w^2 \in \bar{E}^1} \sum_{p^2 \in P_{w^2}^2} x_{p^2}^2 \ln\left(x_{p^2}^2 / d_{w^2}^2\right) \right\},$$

……………………………………………………………………………………..

$$\Phi^k\left(d^k\right) = \min_{\substack{f^k = \Theta^k x^k,\ x^k \in X^k\left(d^k\right) \\ d_{e^k}^{k+1} = f_{e^k}^k,\ e^k \in \bar{E}^k}} \left\{ \sum_{e^k \in \tilde{E}^k} \sigma_{e^k}^k\left(f_{e^k}^k\right) + \Phi^{k+1}\left(d^{k+1}\right) + \gamma^k \sum_{w^k \in \bar{E}^{k-1}} \sum_{p^k \in P_{w^k}^k} x_{p^k}^k \ln\left(x_{p^k}^k / d_{w^k}^k\right) \right\},$$

……………………………………………………………………………………..

$$\Phi^m\left(d^m\right) = \min_{f^m = \Theta^m x^m,\ x^m \in X^m\left(d^m\right)} \left\{ \sum_{e^m \in E^m} \sigma_{e^m}^m\left(f_{e^m}^m\right) + \gamma^m \sum_{w^m \in \bar{E}^{m-1}} \sum_{p^m \in P_{w^m}^m} x_{p^m}^m \ln\left(x_{p^m}^m / d_{w^m}^m\right) \right\}.$$

То, что эти задачи выпуклые, сразу может быть не очевидно. Чтобы это понять, заметим, что ограничения $x^k \in X^k\left(d^k\right)$ с помощью метода множителей Лагранжа можно убрать, добавив в функционал слагаемые

$$\sum_{w^k \in \bar{E}^{k-1}} \max_{\lambda_{w^k}^k} \left\langle \lambda_{w^k}^k, \sum_{p^k \in P_{w^k}^k} x_{p^k}^k - d_{w^k}^k \right\rangle,\ k = 1, \ldots, m.$$

Каждое такое слагаемое – есть выпуклая функция по совокупности параметров $x^k$, $d^k$ (см., например, формулу (3.1.8) стр. 96 [17]). Следовательно (см., например, теорему 3.1.2 стр. 92 и формулу (3.1.9) стр. 96 [17]), $\Phi^m\left(d^m\right)$ – выпуклая функция, но тогда и $\Phi^k\left(d^k\right)$ – выпуклая функция, поскольку (по индукции) $\Phi^{k+1}\left(d^{k+1}\right)$ выпуклая функция ($k = 1, \ldots, m-1$).

**Теорема 1.** 1. *Задачи (2) и (3) являются эквивалентными задачами выпуклой оптимизации, имеющими единственное решение.*

2. *Введенная марковская логит динамика при* $N \to \infty$ *– эргодическая. Ее финальное распределение (возникающее в пределе* $T \to \infty$*) совпадает со стационарным. В предположении (1) стационарное распределение экспоненциально сконцентрировано в окрестности решения задачи (3) (в пределе* $M \to \infty$ *стационарное распределение полностью сосредотачивается на решении задачи (3)).*

3. *Введенная марковская логит динамика при пределах* $N \to \infty$, $M \to \infty$ *описывается СОДУ. В предположении (1) любая допустимая траектория СОДУ (соответствующая вектору корреспонденций* $d^1$*) сходится при* $T \to \infty$ *к решению задачи (3).*

**Замечание 1.** Утверждения 1, 2 теоремы 1 (кроме единственности решения) остаются верными и в предположении, что по части параметров $\gamma^k$ сделаны предельные переходы (от стохастических равновесий к равновесиям Нэша) $\gamma^k \to 0+$ (важно, что эти переходы осуществляются после предельных переходов, указанных в соответствующих пунктах теоремы 1). К этому же результату (с точки зрения того, к какой задаче оптимизации в итоге сводится поиск равновесий) приводит рассмотрение на соответствующих уровнях вместо логит динамик имитационных логит динамик [6, 18].

**Замечание 2.** Утверждения теоремы 1 и замечания 1 остаются верными, если на части ребер (любого уровня) сделать предельные переходы (важно, что эти переходы осуществляются после предельных переходов, указанных в соответствующих пунктах теоремы) вида (предел стабильной динамики [2, 12])

$$\tau_e^\mu(f_e) \xrightarrow[\mu\to 0+]{} \begin{cases} \bar{t}_e, & f_e < \bar{f}_e \\ [\bar{t}_e, \infty), & f_e = \bar{f}_e \end{cases},$$

$$d\tau_e^\mu(f_e)/df_e \xrightarrow[\mu\to 0+]{} 0,\ 0 \le f_e < \bar{f}_e,$$

с дополнительной оговоркой, что существует такой $x \in X$, что условие $f = \Theta x$ совместно с $\{f_e < \bar{f}_e\}_e$. При этом

$$\sigma_e(f_e) = \lim_{\mu\to 0+} \int_0^{f_e} \tau_e^\mu(z)\,dz = \begin{cases} f_e \bar{t}_e, & f_e \le \bar{f}_e \\ \infty, & f_e > \bar{f}_e \end{cases}.$$

Величину $t_e = \lim_{\mu\to 0+} \tau_e^\mu\left(f_e^\mu\right) \ge \bar{t}_e$ можно понимать как затраты на проезд по ребру $e$ (см. также п. 3), а $\lim_{\mu\to 0+} \tau_e^\mu\left(f_e^\mu\right) - \bar{t}_e$ – как дополнительные затраты, приобретенные из-за наличия “пробки” на ребре $e$ [2, 19], возникшей из-за функционирования ребра на пределе пропускной способности $\bar{f}_e$. Эти дополнительные затраты в точности совпадают с множителем Лагранжа к ограничению $f_e \le \bar{f}_e$ [2, 19]. Их также можно понимать как оптимальные платы за проезд (для обычных ребер эти платы равны $f_e\, d\tau_e^\mu(f_e)/df_e$ [3, 20, 21]), взимаемые согласно механизму Викри–Кларка–Гроуса [21].

Если для некоторых $1 \le p \le q \le m$ имеют место равенства $\gamma^p = \ldots = \gamma^q$, то можно свернуть $\Gamma^p$, …, $\Gamma^q$ в один граф $\coprod_{k=p}^{q} \Gamma^k$. Это следует из свойств энтропии (см. свойство 3 § 4 главы 2 [22]).

Далее мы отдельно рассмотрим специальный случай $\gamma^1 = \ldots = \gamma^m = \gamma$. В этом случае мы имеем граф

$$\Gamma = \coprod_{k=1}^{m} \Gamma^k = \left\langle V, E = \coprod_{k=1}^{m} \tilde{E}^k \right\rangle,$$

который имеет всего один уровень, а задача (1) может быть переписана следующим образом

$$\Psi(x)=\sum_{e\in E}\sigma_e(f_e)+\gamma\sum_{w^1\in OD^1}\sum_{p\in\tilde{P}^1_{w^1}}x_p\ln\left(x_p/d^1_{w^1}\right)\to\min_{f=\Theta x,\,x\in X}. \quad (4)$$

**Теорема 2.** *При* $\gamma^1=...=\gamma^m=\gamma$

1. *задачи (2) и (4) являются эквивалентными задачами выпуклой оптимизации, имеющими единственное решение;*
2. *введенная марковская логит динамика при* $N\to\infty$ *– эргодическая. Ее финальное распределение (возникающее в пределе* $T\to\infty$*) совпадает со стационарным, которое представимо в виде (представление Санова)*

$$\sim\exp\left(-\frac{M}{\gamma}\cdot\left(\Psi(x)+o(1)\right)\right),\ M\gg1.$$

*Как следствие, получаем, что стационарное распределение экспоненциально сконцентрировано в окрестности решения задачи (4) (в пределе* $M\to\infty$ *стационарное распределение полностью сосредотачивается на решении задачи (4));*

3. *введенная марковская логит динамика при пределах* $N\to\infty$, $M\to\infty$ *описывается СОДУ. Функция* $\Psi(x)$ *является функцией Ляпунова этой СОДУ (принцип Больцмана). То есть убывает на траекториях СОДУ. Как следствие, любая допустимая траектория СОДУ (соответствующая вектору корреспонденций* $d^1$*) сходится при* $T\to\infty$ *к решению задачи (4).*

**Замечание 3.** К теореме 2 можно сделать замечания аналогичные замечаниям 1, 2 к теореме 1.

## 3. Двойственная задача

Рассмотрим граф

$$\Gamma=\coprod_{k=1}^{m}\Gamma^k=\left\langle V,E=\coprod_{k=1}^{m}\tilde{E}^k\right\rangle.$$

Обозначим через $t_e=\tau_e(f_e)$ (здесь специально упрощаем обозначения, поскольку в виду предыдущего раздела контекст должен восстанавливаться однозначным образом). Запишем в пространстве $t=\{t_e\}_{e\in E}$ *двойственную задачу* к (3) [1, 2, 8] (далее мы используем обозначение $\operatorname{dom}\sigma^*$ – область определения сопряженной к $\sigma$ функции)

$$\min_{f,x}\left\{\Psi(x,f):\ f=\Theta x,\ x\in X\right\}=$$

$$=-\min_{t\in\operatorname{dom}\sigma^*}\left\{\gamma^1\psi^1\left(t/\gamma^1\right)+\sum_{e\in E}\sigma^*_e(t_e)\right\}, \quad (5)$$

где

$$\sigma^*_e(t_e)=\max_{f_e}\left\{f_et_e-\int_0^{f_e}\tau_e(z)dz\right\},$$

$$\frac{d\sigma^*_e(t_e)}{dt_e}=\frac{d}{dt_e}\max_{f_e}\left\{f_et_e-\int_0^{f_e}\tau_e(z)dz\right\}=f_e:\ t_e=\tau_e(f_e),\ e\in E;$$

$$g_{p^m}^m(t)=\sum_{e^m\in\tilde{E}^m}\delta_{e^m p^m}t_{e^m}=\sum_{e^m\in E^m}\delta_{e^m p^m}t_{e^m},$$

$$g_{p^k}^k(t)=\sum_{e^k\in\tilde{E}^k}\delta_{e^k p^k}t_{e^k}-\sum_{e^k\in\bar{E}^k}\delta_{e^k p^k}\gamma^{k+1}\psi_{e^k}^{k+1}\left(t/\gamma^{k+1}\right),\ k=1,\dots,m-1,$$

$$\psi_{w^k}^k(t)=\ln\left(\sum_{p^k\in P_{w^k}^k}\exp\left(-g_{p^k}^k(t)\right)\right),\ k=1,\dots,m,$$

$$\psi^1(t)=\sum_{w^1\in OD^1}d_{w^1}^1\psi_{w^1}^1(t).$$

**Теорема 3.** *Имеет место сильная двойственность (5). Решение задачи выпуклой оптимизации (5) $t\ge 0$ существует и единственно. По этому решению однозначно можно восстановить решение исходной задачи (3) (если какой-то из $\gamma^k\to 0+$, то однозначность восстановления $x$ может потеряться)*

$$f=\Theta x=-\nabla\psi^1\left(t/\gamma^1\right),$$

$$x_{p^k}^k=d_{w^k}^k\frac{\exp\left(-g_{p^k}^k(t)/\gamma^k\right)}{\sum_{\tilde{p}^k\in P_{w^k}^k}\exp\left(-g_{\tilde{p}^k}^k(t)/\gamma^k\right)},\ p^k\in P_{w^k}^k,\ w^k\in OD^k,\ k=1,\dots,m. \tag{6}$$

*Верен и обратный результат. Пусть $f=\Theta x$ – решение задачи (3), тогда $t=\left\{\tau_e\left(f_e\right)\right\}_{e\in E}$ – единственное решение задачи (3) (если какой-то из $\gamma^k\to 0+$, то решение $x$ может быть не единственно, однако это никак не сказывается на возможности однозначного восстановления $t$).*

**Замечание 4.** К теореме 3 можно сделать замечания аналогичные замечаниям 1, 2 к теореме 1. При этом оговорки, возникающие при $\gamma^k\to 0+$, частично уже были сделаны в формулировке самой теоремы. Дополним их следующим наблюдением. Слагаемое $\gamma^1\psi^1\left(t/\gamma^1\right)$ в двойственной задаче (5) имеет равномерно ограниченную константу Липшица градиента в 2-норме:

$$L_2\le\frac{1}{\min_{k=1,\dots,m}\gamma^k}\sum_{w^1\in OD^1}d_{w^1}^1\max_{\breve{p}_{w^1}^1\in\breve{P}_{w^1}^1}\left|\breve{p}_{w^1}^1\right|^2,$$

где $\left|\breve{p}_{w^1}^1\right|$ – число ребер в пути $\breve{p}_{w^1}^1$. Эта гладкость теряется при $\gamma^k\to 0+$:

$$-\lim_{\gamma^k\to 0+}\gamma^k\psi_{w^k}^k\left(t/\gamma^k\right)=\min_{p^k\in P_{w^k}^k}g_{p^k}^k(t)$$

– длина кратчайшего пути в графе $\Gamma^k$, отвечающего корреспонденции $w^k\in OD^k$, ребра $e^k\in\bar{E}^k$, которого взвешены величинами $\gamma^{k+1}\psi_{e^k}^{k+1}\left(t/\gamma^{k+1}\right)$, которые можно понимать как “средние” затраты на $e^k\in\bar{E}^k$ (см. замечание 5). Заметим также, что в пределе (стабильной динамики) $\mu\to 0+$ (см. замечание 2) получаем:

$$\sigma_e^*(t_e)=\lim_{\mu\to 0+}\max_{f_e}\left\{f_e t_e-\int_0^{f_e}\tau_e^{\mu}(z)\,dz\right\}=\begin{cases}\bar{f}_e\cdot\left(t_e-\bar{t}_e\right), & t_e\ge\bar{t}_e\\ \infty, & t_e<\bar{t}_e\end{cases}.$$

При этом $\bar{f}_e-f_e$ в точности совпадает с множителем Лагранжа к ограничению $t_e\ge\bar{t}_e$ [2, 19].

**Замечание 5.** Формулу (6) можно получить и из других соображений. Предположим, что каждый пользователь $l$ транспортной сети, использующий корреспонденцию $w^k \in OD^k$ на уровне $k$ (ребро $e^{k-1}\left(=w^k\right)\in \bar{E}^{k-1}$ на уровне $k-1$), выбирает маршрут следования $p^k \in P_{w^k}^k$ на уровне $k$, если

$$p^k = \arg\max_{q^k\in P_{w^k}^k}\left\{-g_{q^k}^k\left(t\right)+\xi_{q^k}^{k,l}\right\},$$

где независимые случайные величины $\xi_{q^k}^{k,l}$, имеют одинаковое двойное экспоненциальное распределение, также называемое распределением Гумбеля [5, 6, 8]:

$$P\left(\xi_{q^k}^{k,l}<\zeta\right)=\exp\left\{-e^{-\zeta/\gamma^k-E}\right\}.$$

Отметим также, что если взять $E\approx 0.5772$ – константа Эйлера, то

$$M\left[\xi_{q^k}^{k,l}\right]=0,\; D\left[\xi_{q^k}^{k,l}\right]=\left(\gamma^k\right)^2\pi^2\big/6.$$

Распределение Гиббса (логит распределение) (6) получается в пределе, когда число агентов на каждой корреспонденции $w^k\in OD^k$, $k=1,...,m$ стремится к бесконечности (случайность исчезает и описание переходит на средние величины). Полезно также в этой связи иметь в виду, что [5, 9]

$$\gamma^k\psi_{w^k}^k\left(t/\gamma^k\right)=M_{\left\{\xi_{p^k}^k\right\}_{p^k\in P_{w^k}^k}}\left[\max_{p^k\in P_{w^k}^k}\left\{-g_{p^k}^k\left(t\right)+\xi_{p^k}^k\right\}\right].$$

Таким образом, если каждый пользователь сориентирован на вектор затрат $t$ на ребрах $E$ (одинаковый для всех пользователей) и на каждом уровне (принятия решения) пытается выбрать кратчайший путь исходя из зашумленной информации и исходя из усреднения деталей более высоких уровней (такое усреднение можно обосновывать, если, например, как в п. 2, ввести разный масштаб времени (частот принятия решений) на разных уровнях, а можно просто постулировать, что пользователь так действует, как это принято в моделях типа Nested Logit [5]), то такое поведение пользователей (в пределе, когда их число стремится к бесконечности) приводит к описанию распределения пользователей по путям/ребрам (6). Равновесная конфигурация характеризуется тем, что вектор $t$ породил согласно формуле (6) такой вектор $f$, что имеет место соотношение $t=\left\{\tau_e\left(f_e\right)\right\}_{e\in E}$. Поиск такого $t$ (неподвижной точки) приводит к задаче (5).

**Замечание 6.** Сопоставить формуле (4), теореме 2 и замечанию 3 (отвечающих случаю $\gamma^1=...=\gamma^m=\gamma$) вариант двойственной задачи (5) чрезвычайно просто (мы здесь опускаем соответствующие выкладки). Собственно, понять формулу (4) как раз проще не из свойств энтропии (как это было описано в п. 2), а с помощью обратного перехода от двойственного описания (5). Теорема 3 и замечание 5 в случае $\gamma^1=...=\gamma^m=\gamma$ наглядно демонстрируют отсутствие какой бы то ни было иерархии, и возможность работать на одном графе с естественной интерпретацией функций затрат на путях $g_p\left(t\right)$ (без всяких "средних" оговорок).

Перейдем к конспективному обсуждению численных аспектов решения задачи (5). Как правило, выгоднее решать именно задачу (5), а не (3) [8]. На эту задачу удобно смотреть, как на гладкую (с Липшицевым градиентом) задачу композитной оптимизации [9, 10] с евклидовой прокс-структурой (задаваемой 2-нормой). При этом, даже если по ряду параметров $\gamma^k$ требуется сделать предельный переход $\gamma^k \to 0+$, то, как правило, лучше считать, что численно мы все равно решаем задачу со всеми $\gamma^k > 0$ [8]. Этого можно добиться обратным процессом: энтропийной регуляризацией прямой задачи = сглаживанием двойственной. Некоторые детали того, как именно и в каких случаях полезно сглаживать задачу (5) описаны в работе [8] (см. также [23]).

Композитный быстрый градиентный метод (и различные его вариации с адаптивным подбором константы Липшица градиента, универсальный метод и др. [9, 10, 24]) обладает прямо-двойственной структурой [2, 8, 25 – 28]. Это означает, что генерируемые этим методом последовательности $\{t^i\}$ и $\{\tilde{t}^i\}$ обладают следующим свойством

$$\gamma^1\psi^1\left(\tilde{t}^N/\gamma^1\right)+\sum_{e\in E}\sigma_e^*\left(\tilde{t}_e^N\right)-$$

$$-\min_{t\in\operatorname{dom}\sigma^*}\left\{\frac{1}{A_N}\left[\sum_{i=0}^N a_i\cdot\left(\gamma^1\psi^1\left(t^i/\gamma^1\right)+\left\langle\nabla\psi^1\left(t^i/\gamma^1\right),t-t^i\right\rangle\right)\right]+\sum_{e\in E}\sigma_e^*\left(t_e\right)\right\}\le\frac{CL_2R_2^2}{A_N}, \qquad (7)$$

где константа $C \le 10$ зависит от метода,

$$a_N \sim N,\ A_N=\sum_{i=0}^N a_i,\ A_N \sim N^2,$$

$$R_2^2=\max\left\{\tilde{R}_2^2,\hat{R}_2^2\right\},\ \tilde{R}_2^2=\frac{1}{2}\left\|\bar{t}-t^*\right\|_2^2,\ \hat{R}_2^2=\frac{1}{2}\sum_{e\in E}\left(\tau_e\left(\bar{f}_e^N\right)-t_e^*\right)^2,$$

$\bar{f}^N$ определяется в теореме 4, метод стартует с $t^0=\bar{t}$, $t^*$ – решение задачи (5).

**Теорема 4.** *Пусть задача (5) решается прямо-двойственным методом, генерирующим последовательности* $\{t^i\}$ *и* $\{\tilde{t}^i\}$, *с оценкой скорости сходимости (7), тогда*

$$0\le\left\{\gamma^1\psi^1\left(\tilde{t}^N/\gamma^1\right)+\sum_{e\in E}\sigma_e^*\left(\tilde{t}_e^N\right)\right\}+\Psi\left(\bar{x}^N,\bar{f}^N\right)\le\frac{CL_2R_2^2}{A_N},$$

*где*

$$f^i=\Theta x^i=-\nabla\psi^1\left(t^i/\gamma^1\right),\ x^i=\left\{x_{p^k}^{k,i}\right\}_{p^k\in P_{w^k}^k,\,w^k\in OD^k}^{k=1,\dots,m},$$

$$x_{p^k}^{k,i}=d_{w^k}^k\frac{\exp\left(-g_{p^k}^k\left(t^i\right)/\gamma^k\right)}{\sum_{\tilde{p}^k\in P_{w^k}^k}\exp\left(-g_{\tilde{p}^k}^k\left(t^i\right)/\gamma^k\right)},\ p^k\in P_{w^k}^k,\ w^k\in OD^k,\ k=1,\dots,m,$$

$$\bar{f}^N=\frac{1}{A_N}\sum_{i=0}^N a_i f^i,\ \bar{x}^N=\frac{1}{A_N}\sum_{i=0}^N a_i x^i.$$

**Замечание 7.** В общем случае описанный выше подход, представляется наиболее предпочтительным. Однако для различных специальных случаев приведенные оценки, по-видимому, можно немного улучшить [4, 28].

Приведенная теорема 4 оценивает число необходимых итераций. Но на каждой итерации необходимо считать $\nabla\psi^1\left(t/\gamma^1\right)$, а для ряда методов и $\psi^1\left(t/\gamma^1\right)$ (например, для всех адаптивных методов, настраивающихся на параметры гладкости задачи [9, 10, 24]). Подобно [1, 7, 8] можно показать (с помощью сглаженного варианта метода Форда–Беллмана), что для этого достаточно сделать $O\left(\left|O^1\right|\left|E\right|\max_{\breve{p}^1\in\breve{P}^1}\left|\breve{p}^1\right|\right)$ арифметических операций. Однако необходимо обговорить один нюанс. Для возможности использовать сглаженный вариант метода Форда–Беллмана [1, 7, 8] необходимо предположить, что любые движения по ребрам графа с учетом их ориентации являются допустимыми, т.е. множество путей, соединяющих заданные две вершины (источник и сток), – это множество всевозможных способов добраться из источника в сток по ребрам имеющегося графа с учетом их ориентации. Сделанная оговорка не сильно обременительная, поскольку нужного свойства всегда можно добиться раздутием исходного графа в несколько раз за счет введения дополнительных вершин и ребер.

В целом, хотелось бы отметить, что прием, связанный с искусственным раздутием исходного графа путем добавления новых вершин, ребер, источников, стоков является весьма полезным для ряда приложений [2 – 4, 12]. В частности, достаточно популярным является введение фиктивных (с нулевыми затратами) путей-ребер, которые дают возможность ничего не делать пользователям (не перемещаться [3], не торговать [12] и т.п.), что, в свою очередь, позволяет рассматривать ситуации с нефиксированными корреспонденциями $d^1$ [3, 12]. Также популярным приемом является перенесение затрат на преодоление вершин (узлов) графа (перекрестков [3], сортировочных станций [12]) в затраты на прохождение дополнительных ребер, появившиеся при “распутывании” узлов. Но, пожалуй, наиболее важным для большинства приложений является введение фиктивного общего источника и общего стока, соединенных дополнительными ребрами с уже имеющимися вершинами графа [3, 4, 12].

Исследование выполнено при поддержке грантов РФФИ 15-31-70001 мол_а_мос, 15-31-20571-мол_а_вед, 14-01-00722-а, 13-01-12007-офи_м.

Исследование первого автора в части 3 выполнено в ИППИ РАН за счет гранта Российского научного фонда (проект №14-50-00150).

Данная статья представляет собой запись нескольких лекций, прочитанных первым автором в рамках курса “Математическое моделирование транспортных потоков” https://mipt.ru/dcam/science/seminars/matematicheskoe-modelirovanie-transportnykh-potokov-2015.php студентам МФТИ и БФУ им. Канта в осеннем семестре 2015/2016 учебного года.

## Литература